\documentclass[a4, 12pt, leqno]{article}
\usepackage{amsmath}
\usepackage{amssymb}
\usepackage{graphicx}
\usepackage{overpic}
\usepackage{authblk}
\usepackage{setspace}
\usepackage{here}
\usepackage{authblk}
\newtheorem{df}{Definition}[section]

\newtheorem{Th}[df]{Theorem}
\newtheorem{Rm}[df]{Remark}
\newtheorem{lm}[df]{Lemma}

\newcommand{\artanh}{\mathrm{artanh \hspace*{0.6mm}}}

\begin{document}

\title{\bf Head-On Collision of a Pair of Coaxial Circular Vortex Filament}
\author{Masashi A{\sc iki}}
\date{}
\maketitle
\vspace*{-0.5cm}

\begin{abstract}
We consider the head-on collision of two coaxial vortex rings described as the motion of 
two circular vortex filaments under the localized induction approximation. We prove the 
existence of solutions to a system of nonlinear partial differential equations proposed by the 
author in \cite{40} which exhibit head--on collision. We also give a necessary and sufficient 
condition for the initial configuration and parameters of the filaments for head-on collision to occur.
Our results suggest that there exists a critical value \( \gamma_{\ast}>1\)
for the ratio \( \gamma \) of the absolute value of the circulations such that
when \( \gamma\in [1,\gamma_{\ast}]\), two approaching rings will collide, and
when \( \gamma \in (\gamma_{\ast},\infty) \), the ring with the larger circulation
passes through the other and then separate indefinitely.

\end{abstract}

\section{Introduction}
In this paper, we are interested in the head-on collision of two vortex rings sharing the same 
axis of symmetry (coaxial vortex rings) in a incompressible and inviscid fluid. 
For the purposes of this paper, we make a distinction between the two terms
``vortex ring'' and ``circular vortex filament'' by the following.
The term ``vortex ring'' will be used to describe a torus shaped structure in which 
the vorticity of the fluid is concentrated. The term ``circular vortex filament'' will be used 
to describe a circular curve in space for which the vorticity of the fluid is concentrated.
For a vortex filament, the vorticity of the fluid at each point of the curve is directed at the 
direction of the tangent vector. 
Hence, under our terminology, a circular vortex filament can be considered as a 
approximation of a vortex ring in which the core size is taken to be zero.

The study of the interaction of coaxial vortex rings dates back to the pioneering paper by
Helmholtz \cite{41}. In \cite{41}, Helmholtz considered vortex motion in a incompressible
and inviscid fluid based on the Euler equations. 
His study includes the motion of circular vortex filaments, and he observed that 
motion patterns such as head-on collision may occur.
Since then, many researches have been done on head-on collision of coaxial vortex rings, and 
interaction of coaxial vortex rings in general.

Dyson \cite{42,43} further studied the interaction of coaxial vortex rings and proposed a 
system of ordinary differential equations describing the rings. 
From here, we will refer to this model as the Dyson model. Dyson numerically considered the 
head-on collision of two identical rings approaching each other and observed the dynamics of the rings 
as the distance between the two rings decreased. 
In Shariff, Leonard, and Ferziger \cite{51} and Shariff, Leonard, Zabusky, and Ferziger \cite{45},
they extend the method of contour dynamics,
introduced in Zabusky, Hughes, and Roberts \cite{52}, and numerically investigated
the head-on collision of coaxial vortex rings. They also compare and contrast with 
dynamics described by the Dyson model.
Direct numerical simulation of the Navier--Stokes equations was done by 
Stanaway, Shariff, and Hussain \cite{48}. In particular, they considered 
the head-on collision of coaxial vortex rings and the effect of 
viscosity was observed. 
Inoue, Hattori, and Sasaki \cite{53} also conducted numerical simulations of the Navier--Stokes 
equations to investigate the head-on collision of coaxial vortex rings when the 
translational velocity of the rings are varied.

Experiments were conducted by many researchers as well.
Oshima \cite{47} conducted an experiment in which she observed the formation of 
multiple small rings after the initial two coaxial rings collided. This occurs due to the 
reconnection of the two rings and her work was the first to observe this phenomenon. 
A more revealing experiment of the reconnection phenomenon and the related instability of the motion 
of vortex rings was given in Lim and Nickels \cite{49}.
Kambe and Minota \cite{44} conducted experiments to observe the acoustic waves radiated by the 
head-on collision of coaxial vortex rings.
Chu et.al. \cite{46} conducted experiments and numerical calculations of the Navier-Stokes equations
to investigate the head-on collision phenomenon and its relation to the change in enstrophy. 

Shariff and Leonard \cite{56} and Maleshko \cite{57} give an in-depth review of the history 
of the research of vortex rings in which many aspects of motion, including head-on collision,
are addressed.

Although the study of head-on collision of coaxial vortex rings have been done for a very long time,
there are significantly fewer research of the head-on collision phenomenon in a 
mathematically rigorous framework.
Giga and Miyakawa \cite{58} and Feng and \( {\rm \check{S}}\)ver\'ak \cite{59} 
proved the well-posedness of the initial value problem for the Navier--Stokes equations with 
initial data given as vortex rings. In principle, these results give a mathematically rigorous
treatment of the interaction of vortex rings, but extracting the dynamics of specific motion 
patterns through this approach seems difficult.
In Borisov, Kilin, and Mamaev \cite{55}, they analyze the Dyson model to determine the possible 
motion patterns of a pair of coaxial vortex rings for a wide range of configurations, but 
their work doesn't include head-on collision. 

Another approach one may take is to consider circular vortex filaments instead of vortex rings.
By simplifying the structure, it is possible that specific motion patterns can be treated in a 
mathematically rigorous framework, and this is the approach we adopt in this paper.
As far as the author knows, the works by Banica and Miot \cite{60} and 
Banica, Faou, and Miot \cite{61,62} are the only mathematically rigorous results considering
the collision of vortex filaments. They considered the motion of nearly parallel vortex filaments 
described by the model system of partial differential equations proposed by 
Klein, Majda, and Damodaran \cite{63}. Since the model  system is derived by assuming 
that the vortex filaments are nearly straight and parallel, it is not suitable for 
describing motions of circular vortex filaments. 
In \cite{40}, the author proposed a system of nonlinear partial differential equations 
describing the 
interaction of vortex filaments. The paper \cite{40} focused on 
deriving a new system describing the interaction of vortex filaments with general shape and 
proving the existence of solutions corresponding to leapfrogging in the case of 
coaxial circular vortex filaments. 
The aim of this paper is to prove the existence 
of solutions to the model system proposed in \cite{40} which correspond to head-on collision of 
coaxial circular vortex filaments. We also give necessary and sufficient conditions on the filament 
configurations and parameters for head-on collision to occur.
This further shows the capabilities of the model to describe vortex filament motion.
The results of this paper will also imply that the time-global solvability of initial value problems for the 
system doesn't hold in general.

The rest of the paper is organized as follows.
In Section 2, we formulate the problem and state our main theorem.
Section 3 is devoted to the proof of the main theorem.
Finally, in Section 4, we give some discussions and concluding remarks.
In particular, we compare our theoretical results with numerical results obtained by
Inoue, Hattori, and Sasaki \cite{53}.

%-----------------------------------------------------------------------------------------------------
%-----------------------------------------------------------------------------------------------------
%-----------------------------------------------------------------------------------------------------

\section{Problem Setting}
\setcounter{equation}{0}
In \cite{40}, the author proposed the following system of nonlinear partial differential equations.
\begin{align}
\left\{
\begin{array}{l}
\displaystyle
\mbox{\mathversion{bold}$X$}_{t}=
\Gamma_{1}\frac{
\mbox{\mathversion{bold}$X$}_{\xi}\times 
\mbox{\mathversion{bold}$X$}_{\xi \xi}}
{
|\mbox{\mathversion{bold}$X$}_{\xi}|^{3}
}
-
\alpha \Gamma_{2}
\frac{\mbox{\mathversion{bold}$Y$}_{\xi}\times
(\mbox{\mathversion{bold}$X$}-\mbox{\mathversion{bold}$Y$})}
{|\mbox{\mathversion{bold}$X$}-\mbox{\mathversion{bold}$Y$}|^{3}},
\\[5mm]
\displaystyle
\mbox{\mathversion{bold}$Y$}_{t}=
\Gamma_{2}\frac{
\mbox{\mathversion{bold}$Y$}_{\xi}\times 
\mbox{\mathversion{bold}$Y$}_{\xi \xi}}
{
|\mbox{\mathversion{bold}$Y$}_{\xi}|^{3}
}
-
\alpha \Gamma_{1}
\frac{\mbox{\mathversion{bold}$X$}_{\xi}\times
(\mbox{\mathversion{bold}$Y$}-\mbox{\mathversion{bold}$X$})}
{|\mbox{\mathversion{bold}$X$}-\mbox{\mathversion{bold}$Y$}|^{3}}.
\end{array}\right.
\label{model}
\end{align}
Here, \( \mbox{\mathversion{bold}$X$}=\mbox{\mathversion{bold}$X$}(\xi, t)\) and 
\( \mbox{\mathversion{bold}$Y$}=\mbox{\mathversion{bold}$Y$}(\xi ,t) \) are the position vector 
of the filaments parametrized by \( \xi \) at time \( t\), non-zero parameters
\( \Gamma_{1}\) and \(\Gamma_{2}\) are the circulation of the filaments 
\( \mbox{\mathversion{bold}$X$}\) and \( \mbox{\mathversion{bold}$Y$}\) respectively, 
\( \alpha \) is a positive parameter introduced in the derivation of the model,
\( \times \) is the exterior product in the three-dimensional Euclidean space,
and subscripts denote differentiation with the respective variables.
The model system (\ref{model}) was derived from the Biot--Savart law by applying the 
localized induction approximation. The localized induction approximation was
applied to the Biot--Savart law first by Da Rios \cite{64} and later independently by 
Murakami et. al. \cite{65} and Arms and Hama \cite{66} to derive a model equation 
describing the motion of a single vortex filament.
In \cite{40}, the concept of localized induction was applied to the case where two vortex filaments are 
present to derive system (\ref{model}).
We first rescale the time variable by a factor of \( \Gamma_{2}\) and 
arrive at
\begin{align}
\left\{
\begin{array}{l}
\displaystyle
\mbox{\mathversion{bold}$X$}_{t}=
\beta \frac{
\mbox{\mathversion{bold}$X$}_{\xi}\times 
\mbox{\mathversion{bold}$X$}_{\xi \xi}}
{
|\mbox{\mathversion{bold}$X$}_{\xi}|^{3}
}
-
\alpha 
\frac{\mbox{\mathversion{bold}$Y$}_{\xi}\times
(\mbox{\mathversion{bold}$X$}-\mbox{\mathversion{bold}$Y$})}
{|\mbox{\mathversion{bold}$X$}-\mbox{\mathversion{bold}$Y$}|^{3}},
\\[5mm]
\displaystyle
\mbox{\mathversion{bold}$Y$}_{t}=
\frac{
\mbox{\mathversion{bold}$Y$}_{\xi}\times 
\mbox{\mathversion{bold}$Y$}_{\xi \xi}}
{
|\mbox{\mathversion{bold}$Y$}_{\xi}|^{3}
}
-
\alpha \beta 
\frac{\mbox{\mathversion{bold}$X$}_{\xi}\times
(\mbox{\mathversion{bold}$Y$}-\mbox{\mathversion{bold}$X$})}
{|\mbox{\mathversion{bold}$X$}-\mbox{\mathversion{bold}$Y$}|^{3}},
\end{array}\right.
\label{model2}
\end{align}
where \( \beta=\Gamma_{1}/\Gamma_{2}\).

Following \cite{40}, we formulate the problem for a pair of coaxial circular vortex filaments.
Suppose that for some \( R_{1,0},R_{2,0}>0\) and \( z_{1,0},z_{2,0}\in \mathbf{R}\), 
the initial filaments 
\( \mbox{\mathversion{bold}$X$}_{0}\) and \( \mbox{\mathversion{bold}$Y$}_{0}\)
are parametrized by \( \xi \in [0,2\pi )\) as follows.
\begin{align*}
\mbox{\mathversion{bold}$X$}_{0}(\xi )={}^{t}
(R_{1,0}\cos(\xi), R_{1,0}\sin(\xi), z_{1,0}), \quad 
\mbox{\mathversion{bold}$Y$}_{0}(\xi)={}^{t}
(R_{2,0}\cos(\xi), R_{2,0}\sin (\xi), z_{2,0}),
\end{align*}
where we assume that \( (R_{1,0}-R_{2,0})^{2}+(z_{1,0}-z_{2,0})^{2}>0\), 
which means that the two circles are not overlapping. 
Now, we make the ansatz
\begin{align*}
\mbox{\mathversion{bold}$X$}(\xi ,t)={}^{t}
(R_{1}(t)\cos (\xi), R_{1}(t)\sin(\xi), z_{1}(t)), \quad 
\mbox{\mathversion{bold}$Y$}(\xi ,t)={}^{t}
(R_{2}(t)\cos (\xi), R_{2}(t)\sin (\xi), z_{2}(t)),
\end{align*}
for the solution and substitute it into (\ref{model2}). After some calculations, we arrive at the 
following initial value problem for a system of ordinary differential equations.
\begin{align}
\left\{
\begin{array}{l}
\displaystyle
\dot{R_{1}}=-\frac{\alpha R_{2}(z_{1}-z_{2})}{\big( (R_{1}-R_{2})^{2}+(z_{1}-z_{2})^{2}\big)^{3/2}}, \\[7mm]
\displaystyle
\dot{z_{1}}=\frac{\beta}{R_{1}}+\frac{\alpha R_{2}(R_{1}-R_{2})}{\big( (R_{1}-R_{2})^{2}+(z_{1}-z_{2})^{2}\big)^{3/2}}, \\[7mm]
\displaystyle
\dot{R_{2}}=\frac{\alpha \beta R_{1}(z_{1}-z_{2})}{\big( (R_{1}-R_{2})^{2}+(z_{1}-z_{2})^{2}\big)^{3/2}}, \\[7mm]
\displaystyle
\dot{z_{2}}=\frac{1}{R_{2}}-\frac{\alpha \beta R_{1}(R_{1}-R_{2})}
{\big( (R_{1}-R_{2})^{2}+(z_{1}-z_{2})^{2}\big)^{3/2}}, \\[7mm]
(R_{1}(0),z_{1}(0),R_{2}(0),z_{2}(0))=(R_{1,0},z_{1,0},R_{2,0},z_{2,0}).
\end{array}\right.
\label{RZode}
\end{align}
Here, a dot over a variable denotes the derivative with respect to \( t\).
The analysis in \cite{40} shows that head-on collision can only occur when 
the circulation of the filaments have opposite signs, i.e. when \( \beta <0\). 
To simplify the notation, we take \( \gamma =-\beta \) to rewrite 
(\ref{RZode}) to obtain
\begin{align}
\left\{
\begin{array}{l}
\displaystyle
\dot{R_{1}}=-\frac{\alpha R_{2}(z_{1}-z_{2})}{\big( (R_{1}-R_{2})^{2}+(z_{1}-z_{2})^{2}\big)^{3/2}}, \\[7mm]
\displaystyle
\dot{z_{1}}=-\frac{\gamma}{R_{1}}+\frac{\alpha R_{2}(R_{1}-R_{2})}{\big( (R_{1}-R_{2})^{2}+(z_{1}-z_{2})^{2}\big)^{3/2}}, \\[7mm]
\displaystyle
\dot{R_{2}}=-\frac{\alpha \gamma R_{1}(z_{1}-z_{2})}{\big( (R_{1}-R_{2})^{2}+(z_{1}-z_{2})^{2}\big)^{3/2}}, \\[7mm]
\displaystyle
\dot{z_{2}}=\frac{1}{R_{2}}+\frac{\alpha \gamma R_{1}(R_{1}-R_{2})}
{\big( (R_{1}-R_{2})^{2}+(z_{1}-z_{2})^{2}\big)^{3/2}}, \\[7mm]
(R_{1}(0),z_{1}(0),R_{2}(0),z_{2}(0))=(R_{1,0},z_{1,0},R_{2,0},z_{2,0}).
\end{array}\right.
\label{nRZode}
\end{align}
We can further assume without loss of generality that \( \gamma \geq 1 \) since the 
case \( \gamma <1\) is reduced to the case \( \gamma \geq 1 \) by renaming the filaments.
By direct calculation, we see that \( \gamma R_{1}^2 -R_{2}^{2}\) is conserved throughout the motion.
This means that \( (R_{1},R_{2}) \) lies on the set defined by 
\( \gamma R_{1}^2 - R_{2}^{2} =d\), where \( d= \gamma R_{1,0}^2 - R_{2,0}^{2}\). Depending on the value of
\( d\), we can further simplify the system, as well as reduce the number of cases we must 
consider to analyze solutions corresponding to head-on collision of coaxial circular vortex filaments.

When \( d\neq 0\), \( (R_{1},R_{2})\) lies on a hyperbola in the \( R_{1}\)-\( R_{2}\) plane and
the variables \( R_{1}\) and \( R_{2}\) can be 
reduced to one variable. We explain the case \( d>0\) in detail
since the case \( d<0\) is identical. Introducing the change of variables
\begin{align*}
R_{1}(t)=\left(\frac{d}{\gamma}\right)^{1/2}\cosh(\theta(t)), \quad 
R_{2}(t)=d^{1/2}\sinh(\theta(t)), \quad W(t)=z_{1}(t)-z_{2}(t),
\end{align*}
we obtain
\begin{align*}
\left\{
\begin{array}{l}
\displaystyle
\dot{\theta} = 
-\frac{\alpha \gamma^{1/2}W}{
\big( \frac{d}{\gamma}(\cosh \theta -\gamma^{1/2}\sinh \theta)^{2}+W^{2}
 \big)^{3/2}}, \\[7mm]
\displaystyle
\dot{W}=
-\frac{1}{d^{1/2}}
\left(
\frac{\gamma^{3/2}}{\cosh \theta}+\frac{1}{\sinh \theta}
\right)
+
\frac{
\alpha d(\sinh \theta - \gamma^{1/2}\cosh \theta )
(\cosh \theta - \gamma^{1/2}\sinh \theta)}
{\gamma^{1/2}
\big( \frac{d}{\gamma}(\cosh \theta -\gamma^{1/2}\sinh \theta)^{2}+W^{2}
 \big)^{3/2}},
\end{array}\right.
\end{align*}
with appropriate initial data.
The above system is of Hamiltonian form and the Hamiltonian is given by
\begin{align*}
\frac{1}{d^{1/2}}
\bigg(
2\gamma^{3/2}\arctan \big( \tanh(\theta/2)\big)+
\log(\tanh(\theta/2))
\bigg)
+
\frac{\alpha \gamma^{1/2}}
{\big( \frac{d}{\gamma}(\cosh \theta -\gamma^{1/2}\sinh \theta)^{2}+W^{2} \big)^{1/2}}.
\end{align*}
Since the Hamiltonian is divergent at the point \(  (\theta ,W )=(\artanh (1/\gamma^{1/2}),0)  \),
which corresponds to the two filaments colliding, the conservation of the 
Hamiltonian asserts that head-on collision cannot occur in this situation.
Hence, we only need to consider the case \( d=0\).

\

When \( d=0\), we see that \( R_{2}=\gamma^{1/2}R_{1}\), and setting
\begin{align*}
\theta(t):=\log(R_{1}(t)) \quad {\rm and} \quad W(t):=z_{1}(t)-z_{2}(t),
\end{align*}
problem (\ref{nRZode}) reads
\begin{align}
\left\{
\begin{array}{l}
\displaystyle
\dot{\theta} = 
-\frac{\alpha \gamma^{1/2}W}{\big( (\gamma^{1/2}-1)^{2}e^{2\theta}+W^{2}\big)^{3/2}}
=:F_{1}(\theta, W), \\[7mm]
\displaystyle
\dot{W}=-(\gamma + \frac{1}{\gamma^{1/2}})e^{-\theta}+
\frac{\alpha \gamma^{1/2}(\gamma^{1/2}-1)^{2}e^{2\theta}}{\big( (\gamma^{1/2}-1)^{2}e^{2\theta}+W^{2}\big)^{3/2}}
=: F_{2}(\theta ,W),
\end{array}\right.
\label{collide}
\end{align}
with initial data given by \( (\theta(0),W(0))=(\theta_{0},W_{0}) \)
with \( \theta_{0}=\log(R_{1,0})\) and \( W_{0}=z_{1,0}-z_{2,0}\). 
In order to describe the behavior of the two filaments,
it is sufficient to consider the behavior of the solutions to system (\ref{collide}).
From here, we analyze system
(\ref{collide}) as a two-dimensional dynamical system. System (\ref{collide}) is of 
Hamiltonian form
and the Hamiltonian \( {\cal H}(\theta,W)\) is given by
\begin{align}
{\cal H}(\theta,W)=-(\gamma+\frac{1}{\gamma^{1/2}})e^{-\theta} 
+ \frac{\alpha \gamma^{1/2}}{\big( (\gamma^{1/2}-1)^{2}e^{2\theta}+W^{2}\big)^{1/2}}.
\label{Hamiltonian}
\end{align}
The phase space and the possible motion patterns vary depending on the value of \( \gamma \). 
Setting \( H_{0}:={\cal H}(\theta_{0},W_{0}) \), 
the conservation of the Hamiltonian yields
\begin{align*}
H_{0}=-(\gamma+\frac{1}{\gamma^{1/2}})e^{-\theta} 
+ \frac{\alpha \gamma^{1/2}}{\big( (\gamma^{1/2}-1)^{2}e^{2\theta}+W^{2}\big)^{1/2}},
\end{align*}
which can be used to rewrite system (\ref{collide}) as follows.
\begin{align}
\left\{
\begin{array}{l}
\displaystyle
\dot{\theta}=-\frac{1}{\alpha^2 \gamma}
\bigg\{
H_{0}+(\gamma +\frac{1}{\gamma^{1/2}})e^{-\theta }
\bigg\}^{2}
\bigg\{
\alpha^{2}\gamma - (\gamma^{1/2}-1)^{2}e^{2\theta}\big[ H_{0}+(\gamma +\frac{1}{\gamma^{1/2}})e^{-\theta } \big]^{2}
\bigg\}^{1/2}, \\[7mm]
\displaystyle
\dot{W}=-(\gamma +\frac{1}{\gamma^{1/2}})e^{-\theta} + 
\frac{(\gamma^{1/2}-1)^2}{\alpha^{2}\gamma}
\bigg\{
H_{0}+(\gamma +\frac{1}{\gamma^{1/2}})e^{-\theta } 
\bigg\}^{3} e^{2\theta}.
\end{array}\right.
\label{collide2}
\end{align}
In particular, we utilized the relation
\begin{align}
W=
\left\{
\frac{\alpha^{2}\gamma}{\big[ H_{0}+(\gamma + \frac{1}{\gamma ^{1/2}})e^{-\theta}\big]^2}
-(\gamma^{1/2}-1)^{2}e^{2\theta}
\right\}^{1/2}.
\label{hamw}
\end{align}
The above relation and the alternate form of system (\ref{collide}) will be used 
extensively throughout the paper.
In order to state our main theorem, we give two definitions for the types of solution we consider.
\begin{df}
{\rm (Head-on collision and asymmetric collision)}\\
%We define two types of solutions: head-on collision and head-on collision-like solutions.
%
For a finite-time solution \( (\theta ,W)\) of system {\rm (\ref{collide})}, 
we define the following two types of solutions.
In what follows, \( T_{max}\in (0,\infty) \) denotes the maximum existence time of the 
solution in consideration.
\begin{description}
\item[(i)] For \( \gamma =1\), we call a finite-time solution 
\( (\theta, W)\in C^{1}([0,T_{max}))\times C^{1}([0,T_{max})) \)
of system {\rm (\ref{collide})} with initial data \((\theta_{0},W_{0})\) satisfying
\( W_{0}\neq 0 \) a solution corresponding to 
head-on collision if \( W(t)\to 0\) monotonically as
\( t\to T_{max}\).
\item[(ii)] For \( \gamma>1\), we call a finite-time solution 
\( (\theta, W) \in C^{1}([0,T_{max}))\times C^{1}([0,T_{max})) \)
of system {\rm (\ref{collide})} with initial data \((\theta_{0},W_{0})\) satisfying
\( W_{0}\neq 0 \) a solution corresponding to 
asymmetric collision if \( W(t)\to 0\) monotonically as
\( t\to T_{max}\).

\end{description}
We will use the terminology ``colliding solution'' to refer to either or both types of 
solutions, depending on the context.

\end{df}
The behavior of the solution in the above two definitions are the same. The only difference is 
the value of \( \gamma\). We made this distinction because the term ``head-on collision'' seems to 
be mostly used when two identical rings collide. 
The term ``asymmetric collision'' was adopted from \cite{53}.
%

%Also note that the reason we restricted the definition by including \( W_{0}>0\)
%is that because of our choice of the unknown variable \( W =z_{1}-z_{2}\), 
%\( W_{0}>0 \) corresponds to a configuration of the filaments in which they initially
%approach each other along the \( z\)-axis.

\medskip

We state our main theorem.

\begin{Th}
For any \( \alpha \in (0,1) \), there exists a unique \( \gamma_{\ast}\in (1,\infty) \) such that 
the following four statements hold.
\begin{description}
\item[(i)] When \( \gamma=1\), the phase space is \( \mathbf{R}^{2}\setminus (\mathbf{R}\times \{0\}) \)
and the following two statements are equivalent.
\begin{description}
\item[(a)] The solution of {\rm (\ref{collide})} with initial data \( (\theta_{0},W_{0})\in 
\mathbf{R}^{2}\setminus (\mathbf{R}\times \{0\}) \) is a solution corresponding to 
head-on collision.
\item[(b)] \( W_{0}>0\).
\end{description}

\item[(ii)] When \( \gamma \in (1,\gamma_{\ast}) \), the phase space is \( \mathbf{R}^{2}\) and there 
exists \( \theta_{\ast}\in \mathbf{R}\) such that the following 
two statements are equivalent.
\begin{description}
\item[(a)] The solution of {\rm (\ref{collide})} with initial data \( (\theta_{0},W_{0})\in \mathbf{R}^{2}\) corresponds to
asymmetric collision.
\item[(b)] \( W_{0}>0\) and one of the following holds.
\begin{description}
\item[(c)] \( {\cal H}(\theta_{0},W_{0}) \leq 0\).
\item[(d)] \( {\cal H}(\theta_{0},W_{0})>0\) and \( \theta_{0}\leq \theta_{\ast} \).
\end{description}

\end{description}

\item[(iii)] When \( \gamma = \gamma_{\ast}\), the phase space is \( \mathbf{R}^{2}\) and the following two statements are 
equivalent.
\begin{description}
\item[(a)] The solution of {\rm (\ref{collide})} with initial data \( (\theta_{0},W_{0})\in \mathbf{R}^{2}\) corresponds to
asymmetric collision.
\item[(b)] \( W_{0}>0\).
\end{description}
\item[(iv)] When \( \gamma \in (\gamma_{\ast},\infty) \), the phase space is \( \mathbf{R}^{2}\).
In this case, none of the solutions correspond to asymmetric collision, and
all solutions of {\rm (\ref{collide})} with initial data 
\( (\theta_{0},W_{0})\in \mathbf{R}^{2}\) has the following 
properties.
\begin{description}
\item[(a)] \( (\theta, W)\in C^{1}\big( (0,\infty) \big) \times C^{1}\big( (0,\infty) \big) \)
\item[(b)] \( W\) is monotonically decreasing and \( W(t)\to -\infty \) as \( t\to \infty\).

\end{description}
\end{description}

\label{main}
\end{Th}
\begin{Rm} {\rm (Note on the assumption for \( \alpha \) in Theorem \ref{main}) } \\
The parameter \( \alpha \) is introduced when we derived the model 
system {\rm (\ref{model})} in {\rm \cite{40}}. \( \alpha \) is explicitly given by 
\begin{align*}
\alpha = \frac{2\delta }{\log(\frac{L}{\varepsilon } )},
\end{align*}
where \( L>0\) is a cut-off parameter and \( \varepsilon >0\) and \( \delta>0\) are
small parameters introduced in the localized induction approximation.
Hence, although the choice of the upper bound on \( \alpha \) in Theorem {\rm \ref{main}} is 
technical, it is natural to assume that \( \alpha >0 \) is small.

\end{Rm}
Theorem \ref{main} gives a necessary and sufficient condition for head-on collision 
and asymmetric collision to occur.
Note that the solution described in (iv) doesn't correspond to asymmetric collision
since it is a 
time-global solution, even though \( W \) could 
monotonically decrease to zero at some finite time.

%-----------------------------------------------------------------------------------------------------
%-----------------------------------------------------------------------------------------------------
%-----------------------------------------------------------------------------------------------------

\section{Proof of Theorem \ref{main} }
\setcounter{equation}{0}
We first note that since \( F_{1}(\theta, W)\) and \( F_{2}(\theta ,W)\) are smooth with respect to
\( \theta \) and \( W \), the time-local unique solvability of initial value problems
relevant to Theorem \ref{main} is known from general theory of ordinary differential equations. 
We denote the maximum existence time for a solution \( (\theta ,W) \) by \( T_{max}\).

Next, we investigate the equilibria of system (\ref{collide}). This will introduce \( \gamma_{\ast}\)
as stated in Theorem \ref{main}.

\subsection{Equilibria of System (\ref{collide}) }
First we determine the equilibria (or the lack there of) of system (\ref{collide}) with \( \gamma >1\).
In this case, the phase space is \( \mathbf{R}^{2}\) and
from the form of \( F_{1}(\theta ,W)\), we see that any equilibrium must have the form
\( (\theta, 0)\). Hence, we set
\( f(\theta ):=F_{2}(\theta ,0) \)
and investigate the zeroes of \( f\). From direct calculation, we have 
\begin{align*}
f(\theta ) = \frac{e^{-\theta }}{(\gamma^{1/2}-1)}
\bigg\{
-(\gamma +\frac{1}{\gamma^{1/2}})(\gamma^{1/2}-1)+\alpha \gamma^{1/2}
\bigg\}.
\end{align*}
Now we set
\begin{align*}
g(\gamma ):=-(\gamma +\frac{1}{\gamma^{1/2}})(\gamma^{1/2}-1)+\alpha \gamma^{1/2},
\end{align*}
and further setting \( \eta:=\gamma^{1/2}>1 \) we have
\begin{align*}
g(\eta^{2}) = \frac{1}{\eta}
\big( 
-\eta^{4}+\eta^3+\alpha \eta^2-\eta +1
\big).
\end{align*}
After some simple calculus, we see that
for any \( \alpha \in (0,1) \), there exists a unique \( \eta_{\ast}\in (1,\infty) \)
such that \( g(\eta_{\ast}^{2}) =0\). Hence, when \( \gamma=\gamma_{\ast}:=\eta^{2}_{\ast} \),
\( f(\theta )=0 \) for all \( \theta \in \mathbf{R}\) and 
\((\theta ,0)\) is an equilibrium for all \( \theta \in \mathbf{R}\). 
When \( \gamma \in (1,\gamma_{\ast}) \), \( f(\theta )>0\) for all \( \theta\in \mathbf{R}\) and when 
\( \gamma \in (\gamma_{\ast},\infty) \), \( f(\theta )<0\) for all \( \theta \in \mathbf{R}\).
In either cases, there are no equilibria.

When \( \gamma =1\), system (\ref{collide}) reduces to
\begin{align*}
\left\{
\begin{array}{l}
\displaystyle
\dot{\theta}=-\frac{\alpha W}{|W|^{3}} \\[7mm]
\displaystyle
\dot{W}=-2e^{-\theta}
\end{array}\right.
\end{align*}
and the corresponding Hamiltonian reduces to
\begin{align*}
{\cal H}(\theta ,W)=-2e^{-\theta} + \frac{\alpha}{|W|}.
\end{align*}
This implies that the phase space is \( \mathbf{R}^{2}\setminus (\mathbf{R}\times \{0\}) \),
and there is no equilibrium.

\medskip

We summarize the results in the following.
\begin{lm}
For any \( \alpha \in (0,1) \), there exists a unique 
\( \gamma_{\ast}\in (1,\infty) \) such that the following holds.
\begin{description}
\item[(i)] When \( \gamma =1\), the phase space is \( \mathbf{R}^{2}\setminus (\mathbf{R}\times \{0\})\)
and there is no equilibrium.
\item[(ii)] When \( \gamma >1 \), the phase space is \( \mathbf{R}^{2}\) and the following hold.
\begin{description}
\item[(a)] If \( \gamma = \gamma_{\ast} \), \( (\theta ,0) \) for all \( \theta \in \mathbf{R}\) are
equilibria.
\item[(b)] If \( \gamma \neq \gamma_{\ast}\), there is no equilibrium.
\end{description}
\end{description}
\label{lmeq}
\end{lm}
From here, we divide the proof according to the value of \( \gamma \).

%-----------------------------------------------------------------------------------------------------
%-----------------------------------------------------------------------------------------------------
%-----------------------------------------------------------------------------------------------------

\subsection{The Case \( \gamma =1 \)}
In this case, system (\ref{collide}) reduces to 
\begin{align}
\left\{
\begin{array}{l}
\displaystyle
\dot{\theta}=-\frac{\alpha W}{|W|^{3}} \\[7mm]
\displaystyle
\dot{W}=-2e^{-\theta}
\end{array}\right.
\label{g1}
\end{align}
and the Hamiltonian is given by
\begin{align}
{\cal H}(\theta ,W) = -2e^{-\theta} + \frac{\alpha}{|W|}.
\label{hg1}
\end{align}
Particularly, we see that \( W\) is monotonically decreasing.
Hence, for a solution to correspond to head-on collision, \( W_{0}>0\) is necessary.

Conversely, for \( W_{0}>0\), we consider the solution of (\ref{g1}) with initial data
\( (\theta_{0},W_{0})\). Since the Hamiltonian is divergent at \( W=0\), the conservation of the Hamiltonian
implies that \( W(t)>0 \) for all \( t\in [0,T_{max}) \).
Hence, (\ref{g1}) and (\ref{hg1}) is further simplified to 
\begin{align}
\left\{
\begin{array}{l}
\displaystyle
\dot{\theta}=-\frac{\alpha }{W^2} \\[7mm]
\displaystyle
\dot{W}=-2e^{-\theta}
\end{array}\right.
\label{gg1}
\end{align}
and 
\begin{align*}
{\cal H}(\theta ,W) = -2e^{-\theta}+\frac{\alpha }{W}.
\end{align*}
The conservation of the Hamiltonian yields
\begin{align*}
-2e^{-\theta}=H_{0}-\frac{\alpha}{W},
\end{align*}
which can be utilized to decouple system (\ref{gg1}). 
In particular, we have
\begin{align}
\dot{W}=H_{0}-\frac{\alpha }{W}.
\label{w}
\end{align}

When \( H_{0}=0\), equation (\ref{w}) can be explicitly solved to obtain
\begin{align*}
W(t)=\sqrt{W_{0}^{2}-\frac{\alpha }{2}t }.
\end{align*}
This shows that the solution is a finite-time solution with \( \displaystyle T_{max}=\frac{2}{\alpha}W_{0}^{2} \), and 
\( W(t)\to 0 \) monotonically as \( t\to T_{max}\). 
This proves that the solution corresponds to 
head-on collision.

When \( H_{0}\neq 0\), solving equation (\ref{w}) gives the following 
implicit formula for \( W \).
\begin{align*}
\frac{\alpha}{H_{0}^{2}}\log(\alpha -H_{0}W) + \frac{W}{H_{0}}
= t+\frac{\alpha}{H_{0}^{2}}\log(\alpha -H_{0}W_{0})+ \frac{W_{0}}{H_{0}}.
\end{align*}
Setting \( G_{1}(W)=\frac{\alpha}{H_{0}^{2}}\log(\alpha -H_{0}W) + \frac{W}{H_{0}}\),
we have
\begin{align*}
G_{1}(W) = t + G_{1}(W_{0}).
\end{align*}
We see that \( G_{1}(W)\) is monotonically decreasing with respect to \( W\) and 
\( G_{1}(W)\to \frac{\alpha }{H_{0}^{2}}\log \alpha <0 \) as
\( W\to 0 \). This implies that \( G_{1}(W_{0})<\frac{\alpha }{H_{0}^{2}}\log \alpha <0 \).
Hence, \( T_{max} = \frac{\alpha }{H_{0}^{2}}\log \alpha-G_{1}(W_{0}) \) and
\( W(t)\to 0 \) monotonically as \( t\to T_{max}\), and corresponds to head-on collision.

In either cases, the solution corresponds to head-on collision and this proves
(i) of Theorem \ref{main}.

\subsection{The Case \( \gamma \in (1,\gamma_{\ast})\) }
By Lemma \ref{lmeq}, the phase space is \( \mathbf{R}^{2}\).
First we make a few observations.
We see that 
\begin{align*}
{\cal H}(\theta ,0)= e^{-\theta}
\bigg\{ 
-(\gamma +\frac{1}{\gamma^{1/2}})+\frac{\alpha \gamma^{1/2}}{(\gamma^{1/2}-1)}
\bigg\}
=
f(\theta )
\end{align*}
and since \( \gamma \in (1,\gamma_{\ast})\), \( {\cal H}(\theta,0) \) is positive,
monotonically decreasing with respect to \( \theta \), and 
\( {\cal H}(\theta ,0) \to 0 \) as \( \theta \to \infty \).

Now, let \( (\theta ,W) \) be a solution of system (\ref{collide2}) with 
initial data \( (\theta_{0},W_{0}) \) satisfying \( W_{0}\neq 0 \) which corresponds to asymmetric collision,
i.e. we assume that (ii) (a) of Theorem \ref{main} holds.
From standard theory of ordinary differential equations, a finite-time solution must 
converge to a boundary point of the phase space or diverge within the phase space as
\( t \to T_{max}\).
Since \( W(t)\to 0 \) as \( t\to T_{max} \), this implies that 
\( \theta(t) \) must diverge to \( \infty \) or \( -\infty \) as \( t\to T_{max}\).

Suppose \( W_{0}<0\). If \( H_{0}\leq  0\),
we see from the equation for \( W \) in system (\ref{collide2}) that 
\begin{align*}
\dot{W}
&=
-(\gamma +\frac{1}{\gamma^{1/2}})e^{-\theta} + 
\frac{(\gamma^{1/2}-1)^2}{\alpha^{2}\gamma}
\bigg\{
H_{0}+(\gamma +\frac{1}{\gamma^{1/2}})e^{-\theta } 
\bigg\}^{3} e^{2\theta}\\[7mm]
& \leq 
\frac{(\gamma +\frac{1}{\gamma^{1/2}})}{\alpha^2 \gamma}
e^{-\theta}
\bigg\{
-\alpha ^{2}\gamma 
+
(\gamma^{1/2}-1)^{2}(\gamma +\frac{1}{\gamma^{1/2}})^{2} 
\bigg\} \\[7mm]
&=
\frac{(\gamma +\frac{1}{\gamma^{1/2}})}{\alpha^2 \gamma}
e^{-\theta}
\bigg\{
-\alpha \gamma^{1/2}
+
(\gamma^{1/2}-1)(\gamma +\frac{1}{\gamma^{1/2}})
\bigg\}
\bigg\{
\alpha \gamma^{1/2}
+
(\gamma^{1/2}-1)(\gamma +\frac{1}{\gamma^{1/2}})
\bigg\}\\[7mm]
&<0,
\end{align*}
which contradicts \( W(t)\to 0 \) as \( t\to T_{max}\).
If \( H_{0}>0\), we see that
since \( W_{0}<0\), \( W\to 0 \) monotonically implies that
\( W(t)<0\) for all \( t\in [0,T_{max}) \). From the equation for \( \theta \)
in system (\ref{collide}), we see that \( \dot{\theta}>0\), which in turn shows that 
\( \theta(t)\to \infty \) as \( t\to T_{max}\). This contradicts the conservation 
of the Hamiltonian since
\( {\cal H}(\theta ,W)\to 0 \) as \( (\theta ,W)\to (\infty, 0) \). 
Hence, \( W_{0}>0 \) is necessary.

Now suppose \( W_{0}>0\) and \( H_{0}>0\).
Define \( g(\theta )\) by
\begin{align*}
\dot{W}=
-(\gamma +\frac{1}{\gamma^{1/2}})e^{-\theta} + 
\frac{(\gamma^{1/2}-1)^2}{\alpha^{2}\gamma}
\bigg\{
H_{0}+(\gamma +\frac{1}{\gamma^{1/2}})e^{-\theta } 
\bigg\}^{3} e^{2\theta}=: g(\theta ),
\end{align*}
and we look for the zeroes of \( g(\theta )\) to determine the behavior of \( W \).
Further setting \( v:=e^{-\theta }\), we have
\begin{align*}
g(\log (1/v)) = v^{-2}
\bigg\{
-(\gamma + \frac{1}{\gamma^{1/2}})v^3 + \frac{(\gamma^{1/2}-1)^{2}}{\alpha^2 \gamma}
\big[
H_{0}+(\gamma + \frac{1}{\gamma^{1/2}})v\big]^3
\bigg\}.
\end{align*}
A final change of variable given by \( y=(\gamma + \frac{1}{\gamma^{1/2}})v \) shows that
finding the zeroes of \( g\) can be reduced to finding the zeroes of 
\begin{align*}
h(y)=-(\gamma + \frac{1}{\gamma^{1/2}})^{-2}y^3 + \frac{(\gamma^{1/2}-1)^2}{\alpha^2 \gamma}
(H_{0}+y)^3.
\end{align*}
in \( (0,\infty) \).
After some simple calculus, we see that for \( \gamma \in (1,\gamma_{\ast}) \),
there exists a unique \( y_{\ast}\in (0,\infty) \) such that
\( h(y_{\ast})=0 \), \( h(y)>0 \) for all \( y\in (0,y_{\ast})\), and 
\( h(y)<0\) for all \( y\in (y_{\ast},\infty) \).
This in turn shows that there exists a unique \( \theta_{\ast}\in \mathbf{R}\) such that
\( g(\theta_{\ast})=0\), \( g(\theta)<0\) for all \( \theta \in (-\infty, \theta_{\ast}) \), and
\( g(\theta )>0\) for all \( \theta \in (\theta_{\ast},\infty) \).
Hence, for \( W(t) \) to tend to zero monotonically, \( \theta_{0} \leq \theta_{\ast} \) is 
necessary. This shows that (ii) (b) of Theorem \ref{main} holds.

\medskip

Now we show that (ii) (b) implies (ii) (a).
Let \( W_{0}>0\), and we consider the solution of system (\ref{collide2}) with initial data
\( (\theta_{0},W_{0})\). 

When \( H_{0}=0\), the equation for \( \theta \) becomes
\begin{align*}
\dot{\theta}
&=
-\frac{(\gamma + \frac{1}{\gamma^{1/2}})^3}{\alpha^{2}\gamma}
\bigg\{
\frac{\alpha^2 \gamma}{(\gamma + \frac{1}{\gamma^{1/2}})^2}-(\gamma^{1/2}-1)^{2}
\bigg\}^{1/2}
e^{-2\theta}\\[5mm]
& =:-m_{0}e^{-2\theta},
\end{align*}
where \( m_{0}>0\) for \( \gamma \in (1,\gamma_{\ast})\). This can be solved explicitly to obtain
\begin{align*}
\theta(t)=\log \big( \frac{1}{2}e^{2\theta_{0}}-2m_{0}t \big).
\end{align*}
Hence, \( T_{max}=\frac{1}{4m_{0}}e^{2\theta_{0}} \) and 
\( \theta(t)\to -\infty \) monotonically as \( t\to T_{max}\).
Setting \(H_{0}=0 \) in (\ref{hamw}), we have
\begin{align*}
W = \bigg\{
\frac{\alpha^2 \gamma}{(\gamma + \frac{1}{\gamma^{1/2}})^2}-(\gamma^{1/2}-1)^{2}
\bigg\}^{1/2}e^{\theta},
\end{align*}
which shows that \( W(t)\to 0 \) monotonically as \( t\to T_{max}\) and corresponds to 
asymmetric collision.

When \( H_{0} < 0\), 
we have
\begin{align*}
&\alpha^{2}\gamma - (\gamma^{1/2}-1)^{2}e^{2\theta}\big[ H_{0}+(\gamma +\frac{1}{\gamma^{1/2}})e^{-\theta } \big]^{2}\\[3mm]
&=
\bigg\{
\alpha \gamma^{1/2}-(\gamma^{1/2}-1)e^{\theta}[H_{0}+(\gamma + \frac{1}{\gamma^{1/2}})e^{-\theta}]
\bigg\}
\bigg\{
\alpha \gamma^{1/2}+(\gamma^{1/2}-1)e^{\theta}[H_{0}+(\gamma + \frac{1}{\gamma^{1/2}})e^{-\theta}]
\bigg\} \\[3mm]
& \geq \alpha \gamma^{1/2}
\big[
\alpha \gamma^{1/2}-(\gamma^{1/2}-1)(\gamma+ \frac{1}{\gamma^{1/2}})
\big].
\end{align*}
Here, we used the fact that \( H_{0}+(\gamma +\frac{1}{\gamma^{1/2}})e^{-\theta }>0\),
which follows from the definition of the Hamiltonian.
Applying the above estimate to the equation for \( \theta \) in system (\ref{collide2}) yields
\begin{align*}
\dot{\theta } &\leq 
-\frac{\big\{
\alpha \gamma^{1/2}
\big[
\alpha \gamma^{1/2}-(\gamma^{1/2}-1)(\gamma+ \frac{1}{\gamma^{1/2}})
\big]
\big\}^{1/2}}{\alpha^{2}\gamma}
\bigg\{
H_{0}+(\gamma + \frac{1}{\gamma^{1/2}})e^{-\theta}
\bigg\}^{2}\\[3mm]
&=: -m_{1}\bigg\{
H_{0}+(\gamma + \frac{1}{\gamma^{1/2}})e^{-\theta}
\bigg\}^{2}.
\end{align*}
Note that \( m_{1}>0\). To estimate \( \theta \), we compare it with the solution to the 
following initial value problem.
\begin{align}
\left\{
\begin{array}{l}
\displaystyle
\dot{\phi } = -m_{1}\bigg\{
H_{0}+(\gamma + \frac{1}{\gamma^{1/2}})e^{-\phi }
\bigg\}^{2}, \\[5mm]
\phi(0) = \theta_{0}.
\end{array}\right.
\label{phi}
\end{align}
Further setting \( u=\frac{(\gamma + \frac{1}{\gamma^{1/2}})}{|H_{0}|}e^{-\phi } \),
problem (\ref{phi}) can be solved explicitly to obtain the following implicit formula for \( u\).
\begin{align}
-\frac{1}{H_{0}^{2}}
\left\{
\log(\frac{u}{u-1}) - \frac{1}{u-1}
\right\}
=
-m_{1}t
-\frac{1}{H_{0}^{2}}
\left\{
\log(\frac{u_{0}}{u_{0}-1}) - \frac{1}{u_{0}-1}
\right\}
\label{solu}
\end{align}
where \( u_{0}=\frac{(\gamma + \frac{1}{\gamma^{1/2}})}{|H_{0}|}e^{-\theta_{0} } \).
Since \( H_{0}<0\), the definition of the Hamiltonian
implies that \( u>1 \). From the definition of \( u\) and the
fact that \( \dot{\phi} <0 \), 
 \( u(t)\) is monotonically increasing with respect to \( t\).
The right-hand side of (\ref{solu}) tends to zero as 
\( t\to T_{\ast} \), where 
\begin{align*}
T_{\ast}=
-\frac{1}{m_{1}H_{0}^{2}}
\left\{
\log(\frac{u_{0}}{u_{0}-1}) - \frac{1}{u_{0}-1}
\right\}
>0.
\end{align*}
This shows that \( u(t)\to \infty \) as \( t\to T_{\ast}\), which in turn shows that
\( \phi(t)\to -\infty \) as \( t\to T_{\ast} \).
By the comparison principle, \( \theta(t) \leq \phi(t) \) as long as both functions exists,
which implies that \( \theta \) is a finite-time solution with \( T_{max}\leq T_{\ast} \) and
\( \theta(t)\to -\infty \) as \( t\to T_{max}\).
From (\ref{hamw}), we see that
\( W(t)\to 0 \) monotonically as \(t \to T_{max} \), and the solution corresponds to 
asymmetric collision.

Finally, when \( H_{0}>0\) and \( \theta_{0}\leq \theta_{\ast} \), \( \theta(t) \leq \theta_{\ast} \) for 
\( t\in [0,T_{max}) \) since \( \dot{\theta}<0 \).
This in turn implies \( \dot{W}\leq 0\). 
We further estimate
\begin{align*}
\dot{\theta}&=-\frac{1}{\alpha^2 \gamma}
\bigg\{
H_{0}+(\gamma +\frac{1}{\gamma^{1/2}})e^{-\theta }
\bigg\}^{2}
\bigg\{
\alpha^{2}\gamma - (\gamma^{1/2}-1)^{2}e^{2\theta}\big[ H_{0}+(\gamma +\frac{1}{\gamma^{1/2}})e^{-\theta } \big]^{2}
\bigg\}^{1/2}\\[5mm]
& =
-\frac{1}{\alpha^{2} \gamma}
\bigg\{
H_{0}+(\gamma +\frac{1}{\gamma^{1/2}})e^{-\theta }
\bigg\}^{3/2}\\[3mm]
& \qquad \times \bigg\{
\alpha^{2}\gamma
[H_{0}+(\gamma +\frac{1}{\gamma^{1/2}})e^{-\theta }]
 - (\gamma^{1/2}-1)^{2}e^{2\theta}\big[ H_{0}+(\gamma +\frac{1}{\gamma^{1/2}})e^{-\theta } \big]^{3}
\bigg\}^{1/2}\\[5mm]
& =
-\frac{1}{\alpha \gamma^{1/2}}
\bigg\{
H_{0}+(\gamma +\frac{1}{\gamma^{1/2}})e^{-\theta }
\bigg\}^{3/2}
\bigg\{
H_{0}-\dot{W}
\bigg\}^{1/2} \\[5mm]
& \leq
-\frac{H_{0}^{1/2}}{\alpha^2 \gamma}
\bigg\{
H_{0}+(\gamma +\frac{1}{\gamma^{1/2}})e^{-\theta }
\bigg\}^{3/2}\\[5mm]
& \leq
-\frac{H_{0}^{1/2}}{\alpha \gamma^{1/2}}(\gamma + \frac{1}{\gamma^{1/2}})^{3/2}e^{-3\theta /2}
=: - m_{2}e^{-3\theta /2},
\end{align*}
where we substituted the equation for \( W\) of system (\ref{collide2}) in the 
third equality.
Like before, we compare \( \theta \) with the solution of the following problem.
\begin{align*}
\left\{
\begin{array}{l}
\dot{\phi} = -m_{2}e^{-3\phi /2}, \\[3mm]
\phi(0)=\theta_{0}.
\end{array}\right.
\end{align*}
The solution of the above problem is given by
\begin{align*}
\phi (t)=\frac{2}{3}\log \big( e^{3\theta_{0}/2}-\frac{3}{2}m_{2}t \big),
\end{align*}
which implies that \( T_{max}\leq \frac{2e^{3\theta_{0}/2}}{3m_{2}} \) and
\( \theta(t)\to -\infty \) monotonically as \( t\to T_{max}\).
From (\ref{hamw}), \( W(t)\to 0\) monotonically as \( t\to T_{max}\).

Hence, for all cases, the solution \( (\theta ,W)\) corresponds to asymmetric collision, i.e.
(ii) (a) holds. This proves (ii) of Theorem \ref{main}.

\subsection{The Case \( \gamma = \gamma_{\ast} \)}

In this case, we see that \( {\cal H}(\theta ,0) =0 \) for all \( \theta \in \mathbf{R}\).
Since the Hamiltonian is monotonically decreasing with respect to \( |W| \), 
we have \( H_{0}<0 \) when \( W_{0}\neq 0\). 

First we assume (iii) (a) holds.
From the equation for \( W \) in system (\ref{collide2}) we have
\begin{align*}
\dot{W}&=-(\gamma +\frac{1}{\gamma^{1/2}})e^{-\theta} + 
\frac{(\gamma^{1/2}-1)^2}{\alpha^{2}\gamma}
\bigg\{
H_{0}+(\gamma +\frac{1}{\gamma^{1/2}})e^{-\theta } 
\bigg\}^{3} e^{2\theta} \\[5mm]
&\leq 
\frac{
-e^{-\theta}(\gamma + \frac{1}{\gamma^{1/2}})
}{\alpha ^2 \gamma}
\bigg\{
\alpha ^{2}\gamma - (\gamma^{1/2}-1)^{2}(\gamma+\frac{1}{\gamma^{1/2}})^{2}
\bigg\}\\[5mm]
& =0.
\end{align*}
Hence, \( W_{0}>0 \) is necessary.

Now suppose (iii) (b) of Theorem \ref{main} holds, i.e. \( W_{0}>0 \).
The conservation of \( {\cal H}(\theta ,W) \) implies that \( W(t)>0\) for all \( t\in [0,T_{max}) \).
From here, we proceed similarly as before and obtain.
\begin{align*}
\dot{\theta } 
&\leq
-\frac{(\gamma^{1/2}-1)^{1/2}|H_{0}|^{1/2}}{\alpha^{7/4}\gamma^{3/4}}
e^{\theta /2}
\bigg\{
H_{0}+(\gamma + \frac{1}{\gamma^{1/2}})e^{-\theta}
\bigg\}^{2}\\[5mm]
& =: -m_{3}e^{\theta /2}
\bigg\{
H_{0}+(\gamma + \frac{1}{\gamma^{1/2}})e^{-\theta}
\bigg\}^{2}.
\end{align*}
Again, by comparing \( \theta \) with the solution to
\begin{align*}
\left\{
\begin{array}{l}
\displaystyle
\dot{\phi } = 
 -m_{3}e^{\phi /2}
\bigg\{
H_{0}+(\gamma + \frac{1}{\gamma^{1/2}})e^{-\phi }
\bigg\}^{2}\\[5mm]
\phi (0) = \theta_{0},
\end{array}\right.
\end{align*}
we conclude that \( T_{max}\leq -\frac{2}{m_{3}}g_{1}(v_{0}) \), where
\begin{align*}
g_{1}(v_{0})
=
\frac{1}{4(\gamma + \frac{1}{\gamma^{1/2}})^{1/2}|H_{0}|^{3/2}}
\bigg\{
\log \bigg( \frac{v_{0}+1}{v_{0}-1}\bigg) -\frac{2v_{0}}{v_{0}^2-1}
\bigg\},
\end{align*}
\( v_{0}=\frac{(\gamma + \frac{1}{\gamma^{1/2}})^{1/2}}{|H_{0}|^{1/2}}e^{-\theta_{0}/2} \),
and \( \theta(t)\to -\infty \) monotonically as \( t\to T_{max} \).
Furthermore, (\ref{hamw}) shows that \( W(t)\to 0 \) monotonically as \( t\to T_{max}\),
and hence, the solution corresponds to asymmetric collision.
This proves that (iii) of Theorem \ref{main} holds.

\subsection{The Case \( \gamma \in (\gamma_{\ast},\infty) \) }

Let \( (\theta_{0},W_{0})\in \mathbf{R}^{2} \).
In this case, we have \( {\cal H}(\theta ,0) <0 \) for 
all \( \theta \in \mathbf{R} \), and the monotonicity of 
\( {\cal H} \) with respect to \(|W|\) asserts that \( H_{0} <0 \).
Furthermore, \( {\cal H}(\theta ,0) \) is monotonically increasing 
with respect to \( \theta \) and 
\begin{align*}
\lim_{\theta\to -\infty}{\cal H}(\theta ,0) = -\infty \quad 
\text{and} \quad 
\lim_{\theta \to \infty}{\cal H}(\theta ,0) = 0.
\end{align*}
Hence, there exists a unique \( \tilde{\theta}\in \mathbf{R}\) such that
\( {\cal H}(\tilde{\theta},0) = H_{0} \). The conservation of \( {\cal H}\) further
implies that
\( \underline{\theta} \leq \theta(t) \) for all \( t\in [0,T_{max}) \), where
\( \underline{\theta} = \tilde{\theta}-1 \).
From the explicit form of the Hamiltonian (\ref{Hamiltonian}), we see that
\( {\cal H}(\theta ,W) \to 0 \) as \( \theta \to \infty \) uniformly with respect to 
\( W\). Hence, there exists \( \overline{\theta}\in \mathbf{R}\) such that
\({\cal H}(\overline{\theta},W) >H_{0} \) for all \( W \in \mathbf{R}\)
and \( \underline{\theta}<\overline{\theta} \).
From the conservation of the Hamiltonian,
\( \theta(t) \leq \overline{\theta } \) for all \( t\in [0,T_{max}) \).

Summarizing the above arguments, we have shown that \( \underline{\theta }\leq \theta (t) \leq \overline{\theta} \)
for all \( t\in [0,T_{max})\), in other words, \( \theta \) is bounded.
From the equation for \( W \) in system (\ref{collide}), we have
\begin{align*}
\dot{W} \leq 
\bigg\{
-(\gamma + \frac{1}{\gamma^{1/2}}) +
\frac{\alpha \gamma^{1/2}}{(\gamma^{1/2}-1)}
\bigg\}
e^{-\underline{\theta}}
=
f(\underline{\theta})<0
\end{align*}
where the last inequality follows from \( \gamma\in (\gamma_{\ast},\infty) \).
We also have
\begin{align*}
\dot{W} \geq -(\gamma+\frac{1}{\gamma^{1/2}})e^{-\overline{\theta}}.
\end{align*}
These two estimates yield the a priori estimate of the form
\begin{align}
W_{0}-(\gamma+\frac{1}{\gamma^{1/2}})e^{-\overline{\theta}}t
\leq W(t)
\leq W_{0}-|f(\underline{\theta})|t
\label{apriori}
\end{align}
for any \( t \geq 0\). Hence, \( T_{max}=\infty \) and 
the solution is a time-global solution. 
Furthermore, estimate (\ref{apriori}) shows that
\( W(t)\to -\infty \) as \( t\to \infty \). This proves (iv) of Theorem \ref{main},
and finishes the proof of Theorem \ref{main}.
\hfill \( \Box \).

%---------------------------------------------------------------------------
%---------------------------------------------------------------------------
%---------------------------------------------------------------------------

\section{Discussions and Concluding Remarks}
\setcounter{equation}{0}
In this paper, we proved the existence of solutions to 
system (\ref{model}) which correspond to two coaxial circular vortex filaments 
colliding. This was done by reducing the problem to system (\ref{collide})
and analyzing the behavior of the solution in detail. We make some remarks.

\subsection{On the Threshold \( \gamma_{\ast}\)}

In most of the preceding works related to head-on collision of coaxial vortex rings,
rings having circulations with equal absolute value were considered.
In our formulation, this corresponds to \( \gamma =1 \).
The results of this paper suggests that if the absolute value of the circulation of the two rings
are close enough, the two rings should collide. In our formulation, 
this corresponds to \( \gamma \in [1,\gamma_{\ast}] \).
This agrees with the numerical observation made by Inoue, Hattori, and Sasaki in \cite{53}.

In \cite{53}, they conducted a numerical simulation of the 
Navier--Stokes equations and investigated the head-on collision of coaxial vortex rings
with equal size, but varying initial translational velocity. 
As is pointed out in \cite{53}, Saffman \cite{67} showed that
the translational velocity \( U\) of a vortex ring can be approximated by
\begin{align*}
U=\frac{\Gamma}{4\pi R_{0}}
\bigg[
\log \left(
\frac{8R_{0}}{R_{c}}\right)
-0.558
\bigg],
\end{align*}
where \( R_{0}\) is the radius of the ring, \( R_{c}\) is the radius of the core,
and \( \Gamma \) is the circulation of the ring. Hence, varying the value of \( U\) for the two rings
while keeping the size of the rings, i.e. \( R_{0}\) and \(R_{c}\), the same, 
is essentially equivalent to varying the value of the circulation of the two rings.
In \cite{53}, they show numerical results for \( M_{1}/M_{2} =1.0,1.1, 1.33,\) and \(2.0 \) (in their notation,
\( M_{1}/M_{2} \) corresponds to our \( \gamma \)).
When \(  M_{1}/M_{2} =1.0\) and \( 1.1 \), the two rings exhibit collision, and when 
\(  M_{1}/M_{2} =1.33\) and \( 2.0 \), the rings pass through one another. 

On the other hand, when \( \alpha =0.2 \), the threshold \( \gamma_{\ast}\) in our formulation
is numerically obtained as \( \gamma_{\ast}=1.219 \). Hence,
Theorem \ref{main} can then be interpreted as follows.
When \( \gamma \in [1,\gamma_{\ast}] \), the two filaments collide, and when 
\( \gamma \in (\gamma_{\ast},\infty) \), the two filaments pass through one another.
This agrees with the numerical findings of \cite{53}.

\subsection{On the Equilibria of System (\ref{collide})}

As is stated in Lemma \ref{lmeq}, the phase space for system (\ref{collide}) is 
\( \mathbf{R}^{2}\) and \( (\theta ,0)\) for all \( \theta \in \mathbf{R}\) 
are equilibria when \( \gamma = \gamma_{\ast}\).
These equilibria correspond to two coaxial circular filaments in perfect balance 
such that there is no relative motion. The effect of self-induction and 
interaction exactly balance each other out and the two filaments stay in the 
same plane throughout the motion. Since this type of behavior is only achieved 
when the two filaments start their motion in the same plane and the two filaments
have circulations with a specific ratio \( \gamma_{\ast}\), it seems 
impractical to recreate such behavior in an experimental setting, even though 
it is suggested in theory.

\subsection{Further Application of the Model System (\ref{model})}

In Chen, Wang, Li, and Wang \cite{68} and Cheng, Lou, and Lim \cite{69}, experiments and 
numerical simulations of the collision of elliptic vortex rings are conducted, respectively.
Unlike a circular vortex ring, an elliptic vortex ring travels while changing shape.
This makes it hard to analyze the phenomenon under the Dyson model, or other 
classical models since they are formulated as systems of ordinary differential equations 
under the assumption that the ring is circular. 

Since system (\ref{model}) is a system of partial differential equations
which describe the motion of filaments with general shape, there is potential for 
system (\ref{model}) to be utilized to conduct mathematical analysis on 
head-on collision of elliptic vortex rings. For a pair of elliptic vortex filaments, 
the axisymmetry of the problem is broken and the problem can no longer be reduced to a system of 
ordinary differential equation. 

The author is currently working on proving the 
solvability of initial value problems to system (\ref{model}). This will serve as a
starting point for mathematically analyzing various phenomena related to interaction of vortex filaments,
including head-on collision of elliptic vortex rings.

\vspace*{1cm}
\noindent
Masashi Aiki\\
Department of Mathematics\\
Faculty of Science and Technology, Tokyo University of Science\\
2641 Yamazaki, Noda, Chiba 278-8510, Japan\\
E-mail: aiki\verb|_|masashi\verb|@|ma.noda.tus.ac.jp

\end{document}